\documentstyle[12pt,leqno,amssymb,graphics]{amsart}

\newtheorem{thm}{Theorem}[section]
\newtheorem{lemma}[thm]{Lemma}
\newtheorem{prop}[thm]{Proposition}

\theoremstyle{definition}
\newtheorem{dfn}[thm]{Definition}
\theoremstyle{remark}
\newtheorem{remark}[thm]{Remark}

\begin{document}

\newcommand{\pr}{\protect\ref}
\newcommand{\su}{\subseteq}
\newcommand{\x}{\times}
\newcommand{\pa}{{\partial}}

\newcommand{\R}{{\Bbb R}}
\newcommand{\Z}{{\Bbb Z}}
\newcommand{\E}{{{\Bbb R}^3}}
\newcommand{\C}{{{\Bbb Z}/2}}

\newcommand{\g}{{{\mathrm{genus}}(F)}}
\newcommand{\A}{{{\mathrm{Arf}}(g)}}
\newcommand{\Ai}{{{\mathrm{Arf}}(g^i)}}
\newcommand{\Ar}{{\mathrm{Arf}}}
\newcommand{\r}{{\mathrm{rank}}}
\newcommand{\F}{{\mathbf{F}}}
\newcommand{\cc}{{\mathbf{C}}}
\newcommand{\nc}{{\mathbf{N}}}

\newcommand{\cd}{{\mathrm{codim}}}

\newcommand{\p}{{\psi}}
\newcommand{\hp}{{\Psi}}

\newcommand{\M}{{\mathcal{M}}}
\newcommand{\hM}{{\widehat{\mathcal{M}}}}
\newcommand{\m}{{\mathcal{M}}_g}
\newcommand{\hm}{{\hM_g}}
\newcommand{\hmi}{{\hM_{g^i}}}
\newcommand{\mi}{{\M_{g^i}}}

\newcommand{\hh}{{H_1}}
\newcommand{\ov}{{O(V,g)}}
\newcommand{\oh}{{O(\hh,g)}}
\newcommand{\ohi}{{O(\hh,g^i)}}

\newcommand{\tta}{{T\circ T_a}}
\newcommand{\T}{{{\mathcal{T}}_c}}
\newcommand{\TT}{{\mathcal{T}}}
\newcommand{\U}{{\mathcal{U}}}

\newcommand{\ep}{{\epsilon}}
\newcommand{\wb}{{W^\bot}}
\newcommand{\w}{{\bigcup_i w_i}}

\newcounter{numb}

\title[Mapping Class Group and Quadruple points]
{Subgroups of the Mapping Class Group \\ and
Quadruple points of Regular Homotopies}
\author{Tahl Nowik}
\address{Department of Mathematics, Columbia University, New York, 
NY 10027, USA.} 
\email{tahl@@math.columbia.edu}
\date{September 2, 1999}

\begin{abstract}

Let $F$ be a closed orientable surface. If $i,i':F\to\E$
are two regularly homotopic generic immersions, then it has been shown in 
[N] that all generic regular homotopies between $i$ and $i'$ have the 
same number mod 2 of quadruple points. We denote this number by $Q(i,i')\in\C$.
We show that for any generic immersion $i:F\to\E$ and any diffeomorphism 
$h:F\to F$ such that $i$ and $i\circ h$ are regularly homotopic,
$$Q(i,i\circ h)=\bigg(\r(h_*-Id)+(n+1)\ep(h)\bigg)\bmod{2},$$
where $h_*$ is the map induced by $h$ on $H_1(F,\C)$, $n$ is the genus
of $F$ and $\ep(h)$ is 0 or 1 according to whether $h$ is orientation
preserving or reversing, respectively.

\end{abstract}

\maketitle

\section{Introduction}\label{INT}

For $F$ a closed surface and $i,i':F\to\E$ two 
regularly homotopic generic immersions, we are interested in the number 
mod 2 of quadruple points occurring in generic regular homotopies 
between $i$ and $i'$.
It has been shown in [N] that this number is the same for all such
regular homotopies, and so it is a function of $i$ and 
$i'$ which we denote $Q(i,i')\in\C$. There then arises the problem of finding 
explicit formulae for $Q(i,i')$. 

Assuming $F$ is orientable, we give an
explicit formula for $Q(i,i\circ h)$, where $i:F\to\E$ is any
generic immersion and $h:F\to F$ is any diffeomorphism such that 
$i$ and $i\circ h$ are regularly homotopic (Theorem \pr{main1}).
For two special cases a formula for $Q(i,i\circ h)$ has already been known:
The case where $F$ is a sphere has appeared in [MB] and [N], and the 
case where $F$ is a torus and $i$ is an embedding has appeared in [N].

Based on the Smale-Hirsch Theorem, Pinkall in [P] gave a useful tool for
determining when two immersions are regularly homotopic, namely, any immersion
$i:F\to\E$ induces a quadratic form $g^i:H_1(F,\C)\to\C$, and two
immersions $i,i':F\to\E$ are regularly homotopic iff $g^i=g^{i'}$.
Let $\hM$ denote the group of all diffeomorphisms $h:F\to F$ up to isotopy.
Given $i:F\to\E$ we are interested in the group of all $h\in\hM$
such that $Q(i,i\circ h)$ is defined, 
that is the group of all $h\in\hM$ such that
$i$ and $i\circ h$ are regularly homotopic. It follows from  
the above criterion that this is precisely the group $\hmi$ of all $h\in\hM$ 
which preserve the quadratic form $g^i$ on $H_1(F,\C)$.
We are thus lead to study the groups $\hm$, starting with their index 2 
subgroup $\m$ of orientation preserving maps.

The plan of the paper is as follows:
In Section \pr{A} we present the known results on quadratic forms
which we will need. In Section \pr{B} we show that the expression
$\r(T-Id)\bmod{2}$ appearing in our proposed formula for $Q(i,i\circ h)$ 
defines a homomorphism on appropriate subgroups of $GL(H_1(F,\C))$
(Theorem \pr{t1}).
In Section \pr{C} we show that (except for one special case)
the group $\m$ is generated by Dehn twists and squares of Dehn twists
(Theorem \pr{t3}).
Our main result is then proved in Section \pr{D}. 
We start with surfaces of genus 0 and 1 where we compute $Q(i,i\circ h)$
by giving explicit regular homotopies for generators of $\hmi$.
We then continue by induction on the genus, using the special nature of the
generators of $\mi$ that we have found and of an additional
generator for $\hmi$. Each generator respects a separation 
of $F$ into surfaces of smaller genus, a fact which
will enable us to construct regular homotopies for $F$ by combining
regular homotopies of surfaces of smaller genus.

\section{Quadratic Forms over $\C$}\label{A}

In this section we summarize the definitions and known properties
of quadratic forms over $\C$ which will be needed in our work.
Proofs to all facts stated in this section may be found in [C], 
except for those relating to the Arf invariant, which may be found in [L].

Let $V$ be a finite dimensional vector space over $\C$. 
A function $g:V\to\C$ is called a \emph{quadratic form} if $g$ satisfies:
$g(x+y)=g(x)+g(y)+B(x,y)$ for all $x,y \in V$, 
where $B(x,y)$ is a bilinear form.
The following properties follow:
(a) $g(0)=0$. (b) $B(x,x)=0$ for all $x\in V$.
(c) $B(x,y)=B(y,x)$ for all $x,y\in V$.  
$g$ is called \emph{non-degenerate} if $B$ is non-degenerate,
i.e. for any $0\neq x\in V$ there is $y\in V$ with $B(x,y)\neq 0$.

\begin{prop}\label{pc1}
If $g$ is non-degenerate then $V$ is necessarily of even dimension and 
there exists
a basis $a_1,\dots,a_n,b_1,\dots,b_n$ for $V$ such that 
$B(a_i,a_j)=B(b_i,b_j)=0$ and 
$B(a_i,b_j)=\delta_{ij}$ for all $1\leq i,j \leq n$
and such that one of the following two possibilities holds:
\begin{enumerate}
\item $g(a_i)=g(b_i)=0$ for $i=1 \dots n$.
\item $g(a_1)=g(b_1)=1$ and $g(a_i)=g(b_i)=0$ for $i=2 \dots n$.
\end{enumerate}
\end{prop}

$g$ is completely determined by the values $g(v_i)$ and $B(v_i,v_j)$ 
on a basis $v_1,\dots,v_{2n}$ 
and so for given dimension $2n$ there are two isomorphism 
classes of non-degenerate quadratic forms, and they are in fact distinct.
The invariant $\A\in\C$ is then defined to be 0 or 1 according to whether 
1 or 2 of Proposition \pr{pc1} holds respectively.
(In the more general setting of [C], this is equivalent to $g$ having index 
$n$ or $n-1$ respectively.) The Arf invariant is additive in the following 
sense: 

\begin{prop}\label{pc2}
If $g_i:V_i\to\C$, $i=1,2$, are non-degenerate quadratic forms,
then $g_1\oplus g_2 :V_1\oplus V_2\to\C$ 
defined by $(g_1\oplus g_2) (x_1,x_2)=g_1(x_1)+g_2(x_2)$
is a non-degenerate
quadratic form with $\Ar(g_1\oplus g_2)=\Ar (g_1)+\Ar (g_2)$.
\end{prop}

From now on we 
will always assume that our quadratic form $g$ is non-degenerate.

\begin{prop}\label{pc3}
If $a_1,\dots,a_k\in V$ are independent and $B(a_i,a_j)=0$ for all
$1\leq i,j \leq k$ then there are $b_1,\dots,b_k\in V$ with
$B(b_i,b_j)=0$ and $B(a_i,b_j)=\delta_{ij}$
for all $1\leq i,j \leq k$. 
($a_1,\dots,a_k,b_1,\dots,b_k$ are then necessarily independent.)
\end{prop}

A linear map $T:V\to V$ is called \emph{orthogonal}
with respect to $g$ if $g(T(x))=g(x)$ for all $x\in V$.
It then follows that $B(T(x),T(y))=B(x,y)$ for all $x,y\in V$
and that $T$ is invertible.
The group of all orthogonal maps of $V$ with respect to $g$ will
be denoted $\ov$.

\begin{dfn}\label{d0}
Given $a\in V$, define $T_a:V\to V$ by $T_a(x)=x+B(x,a)a$.
\end{dfn}

\begin{prop}\label{pc4}
$T_a \in \ov$  iff $g(a)=1$ or $a=0$.
\end{prop}

\begin{thm}[Cartan, Dieudonne]\label{tc}
Except for the case when $\dim V =4$ and $\A=0$,
$\ov$ is generated by the elements $T_a$ with $g(a)=1$.
\end{thm}

Theorem \pr{tc} will also follow from
Theorem \pr{t3} below. See Remark \pr{r5}.

If $W\su V$ is a subspace, then the conjugate space $\wb$ of $W$ is defined
by $W^{\bot}=\{x\in V  : B(x,y)=0 \ \hbox{for all} \ y\in W\}$.
If $a\in V$ we similarly define
$a^{\bot}=\{x\in V  : B(x,a)=0 \}$.
Let $Id$ denote the identity map on $V$, let 
$Im(T)$ denote the image of $T$ and let $\F(T)=\{x\in V : T(x)=x\}$.

\begin{prop}\label{pc5}
If $T\in\ov$ then 
$Im(T-Id) = (\F(T))^\bot$.
\end{prop}

\section{A Homomorphism from $\ov$ to $\C$.}\label{B}

Let $a\in V$, then $\F(T_a)=a^\bot$ and so if $a\neq 0$
$\dim\F(T_a)=2n-1$, where $2n=\dim V$.

\begin{lemma}\label{p1}
Let $T\in \ov$ and $a\in V$ with $g(a)=1$.
\begin{enumerate}
\item If $\F(T)\su \F(T_a)$ then
$\dim\F(\tta) = \dim\F(T)+1$. 
\item If $\F(T)\not\su \F(T_a)$ then
$\dim\F(\tta) = \dim\F(T)-1$. 
\end{enumerate}
\end{lemma}

\begin{pf}
We first note:
(a) If $x\not\in a^\bot$ then $B(x,a)=1$ so $T_a(x) =x+a$ 
and so $\tta(x)=T(x+a)$.
(b) For any $x\in V$, if $T(x+a)=x$ 
then $g(x)=g(T(x+a))=g(x+a)=g(x)+g(a)+B(x,a)$ and so $B(x,a)=g(a)=1$,
that is $x\not\in a^\bot$.

We get that $\F(\tta)\not\su a^\bot$
iff $\exists x\not\in a^\bot$ with $T(x+a)=x$
iff $\exists x\in V$ with $T(x+a)=x$ iff 
$\exists x\in V$ with $(T-Id)(x+a)=a$ iff $a\in Im (T-Id)=\F(T)^\bot$
(Proposition \pr{pc5})
iff $\F(T)\su a^\bot$. Since $a^\bot = \F(T_a)$ we conclude that  
$\F(\tta)\not\su \F(T_a)$ iff $\F(T)\su\F(T_a)$.

Clearly $\F(T)\cap\F(T_a)=\F(\tta)\cap\F(T_a)$ and let $k$ denote 
the dimension of this subspace. 
$\F(T_a)$ is of codimension 1 in $V$ and so it follows:
(1) If $\F(T)\su\F(T_a)$ then $\dim\F(T)=k$ and $\F(\tta)\not\su\F(T_a)$
and so $\dim\F(\tta)=k+1$. (2) If $\F(T)\not\su\F(T_a)$ then $\dim\F(T)=k+1$
and $\F(\tta)\su\F(T_a)$ and so $\dim\F(\tta)=k$.
\end{pf}

We now define $\p:\ov\to\C$ by: 
$$\p(T) \ = \ \r(T-Id)\mod{2}.$$

\begin{remark}\label{r}
Since $\F(T)=\ker(T-Id)$ (or by Proposition \pr{pc5}) we may also write:
$\p(T)=\cd\F(T)\bmod{2}$, and since
$V$ is of even dimension we also have:
$\p(T)=\dim\F(T)\bmod{2}$.
\end{remark}

\begin{thm}\label{t1}
$\p:\ov\to\C$ is a (non-trivial) homomorphism.
\end{thm}

\begin{pf}
We will be using the equivalent definition
$\p(T)=\dim\F(T)\bmod{2}$ of Remark \pr{r}.
Assume first that $(V,g)$ is not of the special case excluded from
Theorem \pr{tc}, and so $\ov$ is generated
by the elements $T_a$ with $g(a)=1$. If $T=T_{a_1}\circ\cdots\circ T_{a_k}$
($g(a_i)=1$) then since $\p(Id)=\dim V \bmod{2}=0$, induction on
Lemma \pr{p1} implies $\p(T)=k\bmod{2}$ 
which clearly implies that $\p$ is a homomorphism. 

We are left with the case $\dim V =4, \A=0$.
By Proposition \pr{pc2}, $(V,g)\cong(V'\oplus V',g'\oplus g')$
where $\dim V'=2,\Ar(g')=1$. We identify $V$ with $V'\oplus V'$ via such an
isomorphism. The set of all elements in $V$ with $g=1$ is $V_1\cup V_2$ 
where $V_1=\{(x,0) : 0\neq x \in V'\}$ 
  and $V_2=\{(0,x) : 0\neq x \in V'\}$.
If $a\in V_1$ and $b\in V_2$ then $B(a,b)=0$, 
whereas if $a\neq b$ are in the same $V_k$ then
$B(a,b)=1$. It follows that any $T\in\ov$ must either map each $V_k$ 
into itself or map $V_1$ into $V_2$ and $V_2$ into $V_1$.
So $T$ is of the form 
$(x,y)\mapsto (T_1(x),T_2(y))$
or $(x,y)\mapsto (T_1(y),T_2(x))$ where $T_1,T_2\in O(V', g')$ ($=GL(V')$.)
Such a map will be denoted by $(T_1,T_2)_0$ or $(T_1,T_2)_1$ respectively.
If $T=(T_1,T_2)_0$ then $T(x,y)=(x,y)$ iff $T_1 (x) = x$ and $T_2 (y) = y$
and so $\F(T)=\F(T_1)\oplus\F(T_2)$ so 
$\dim\F(T)=\dim\F(T_1)+\dim\F(T_2)$ so 
$\p(T)=\p(T_1)+\p(T_2)$. 
(The $\p$ on the left is the function on $\ov$ and the $\p$ on the 
right is the function on $O(V',g')$.)
If $T=(T_1,T_2)_1$ then
$T (x,y)=(x,y)$
iff $T_1 (y)=x$ and $T_2 (x)=y$ that is  $(x,y)$ is of the form $(x,T_2(x))$ 
with $T_1 \circ T_2 (x) = x$. And so $\dim\F(T)=\dim\F(T_1\circ T_2)$
so $\p(T)=\p(T_1\circ T_2)$.
Now, since $\dim V' = 2$, $V'$ belongs to the general case,
and so we already know $\p(T_1 \circ T_2)=\p(T_1) + \p(T_2)$.
So we have shown for both $u=0$ and $u=1$ that
$\p((T_1,T_2)_u)=\p(T_1)+\p(T_2)$ . 
Now if $T=(T_1,T_2)_u$ and
$S=(S_1,S_2)_{u'}$ then $T\circ S$ is of the form
$(T_1\circ S_1 , T_2 \circ S_2)_{u''}$ or
$(T_1\circ S_2 , T_2 \circ S_1)_{u''}$. In any case (again using the fact
that $\p$ on $V'$ is a homomorphism) we get that
$\p(T\circ S)=\p(T_1)+\p(T_2)+\p(S_1)+\p(S_2)=\p(T)+\p(S)$.

Finally, $\p:\ov\to\C$ 
is not trivial since $\p(T_a)=1$ for any $a\in V$ with $g(a)=1$.

\end{pf}

\begin{remark}\label{r1}
1. For $A\in O_k(\R)$ (the group of $k\x k$ orthogonal matrices over $\R$) 
$\cd\F(A)=0\bmod{2}$ iff $\det A =1$.
And so by Remark \pr{r},
$\p:\ov\to\C$ may be thought of as an analogue 
of the homomorphism $\det:O_k(\R)\to\{1,-1\}$.
($\det$ on $\ov$ is of course trivial.)

2. Our expression for $\p$ is meaningful on the whole of $GL(V)$,
however $\p$ is in general not a homomorphism on $GL(V)$ or
even on its subgroup $Sp(V) \supseteq \ov$ 
of maps preserving $B(x,y)$.
\end{remark}

For $\dim V = 4,\A=0$,
we note that though the identification of $V$ with $V'\oplus V'$ 
in the proof of Theorem \pr{t1} is not unique,
the (unordered) pair of sets $V_1,V_2$ is uniquely defined by its 
mentioned properties, namely, $V_1\cup V_2 = \{v\in V : g(v)=1\}$, 
$B(a,b)=0$ for $a\in V_1, b\in V_2$, and
$B(a,b)=1$ for $a\neq b \in V_k$, $k=1,2$.
It follows, as we have noticed, that any $T\in\ov$ either preserves
each $V_k$ (then $T=(T_1,T_2)_0$,) or interchanges the $V_k$s 
(then $T=(T_1,T_2)_1$.)

\begin{dfn}\label{u}
Let $\dim V =4,\A=0$. $T\in\ov$ will be called a 
\emph{$U$-map} if $T$ interchanges $V_1$ and $V_2$.
\end{dfn}

\begin{lemma}\label{pum}
Let $\dim V =4,\A=0$. If $T$ is a $U$-map such that $T^2=Id$
then $\p(T)=0$.
\end{lemma}

\begin{pf} 
$T=(T_1,{T_1}^{-1})_1$ and so by the proof of Theorem \pr{t1},
$\p(T)=\p(T_1)+\p({T_1}^{-1})=0$.
\end{pf}

\section{Generators for the Orthogonal Mapping Class Group}\label{C}

Let $F$ be a closed orientable surface. $\hh$ from now on will always
denote $H_1(F,\C)$ (considered as a vector space over $\C$.) 
Let $g:\hh\to\C$ be a quadratic
form whose associated bilinear form $B(x,y)$ is the algebraic
intersection form $x\cdot y$ of $\hh$. (In particular, $g$ is non-degenerate.)
Let $\M$ 
denote the mapping class group of $F$ i.e. the group of all orientation 
preserving diffeomorphisms $h:F\to F$ up to isotopy.
For $h:F\to F$, let $h_*$ denote the map it induces on $\hh$.
The \emph{orthogonal mapping class group of $F$
with respect to $g$} will be the subgroup $\m$ of $\M$ defined by
$\m=\{h\in\M : h_*\in \oh \}$.

A simple closed curve will be called a \emph{circle}.
If $c$ is a circle in $F$, the homology class of $c$ 
in $\hh$ will be denoted
by $[c]$. Given a circle in $F$, a Dehn twist along $c$ will be denoted $\T$.
(We will not establish a convention as to which of the two possible Dehn 
twists is $\T$ and which is $\T^{-1}$, rather, it will be clear in each 
case which of the two should be used.) The map induced on $\hh$ by $\T$
is $T_{[c]}$ of Definition \pr{d0}.
And so by Proposition \pr{pc4},
$\T\in\m$ iff $g([c])=1$ or
$[c]=0$. Also, since $(T_{[c]})^2=Id$, $(\T)^2\in\m$ for any circle $c$.
In view of this we make the following definition:

\begin{dfn}\label{d1}
A map $h:F\to F$ will be called \emph{good} if it is of
one of the following forms:
\begin{enumerate}
\item $h=(\T)^2$ for some circle $c$.
\item $h=\T$ for a circle $c$ with $g([c])=1$. 
\item $h=\T$ for a circle $c$ with $[c]=0$.
\end{enumerate}
A good map will be called of type 1, 2 or 3 accordingly.
\end{dfn}

The purpose of this section is to show that except for the special case when
$\g=2$ and $\A=0$, $\m$ is generated by the good maps.
For the mentioned special case, 
we will show that one more generator is required. 

Whenever we consider two circles in $F$, we will assume that they intersect
transversally. $|c_1\cap c_2|$ will then denote the number of intersection
points between circles $c_1$ and $c_2$. (And so the algebraic
intersection $[c_1]\cdot [c_2]$ 
in $\hh$ is the reduction mod 2 of $|c_1\cap c_2|$.)
Given two circles $a,b$ in $F$ with $|a\cap b|=1$, there are two ways for 
joining them into one circle $c$ via surgery at their intersection point.
$c$ will be called a \emph{merge} of $a$ and $b$.
If $a$ and $b$ are oriented, then $c$ will be called the \emph{positive}
or \emph{negative} merge of $a$ and $b$, according to whether
the surgery is performed so that the orientations of $a$ and $b$ match or
do not match, respectively.

\begin{lemma}\label{l1}
Let $a,b$ be two oriented circles in $F$ with $|a \cap b|=1$,
let $P$ be their intersection point and let $c$ be their negative merge. 
\begin{enumerate} 
\item $\T$ followed by
an isotopy performed in a thin neighborhood of $a\cup b$,
maps $a$ orientation preservingly onto $b$.
\item If $d$ is another circle passing $P$, and otherwise disjoint form
$a$ and $b$, and if $a$ and $b$ cross $d$ at $P$ in the same direction, 
then the above Dehn twist and isotopy may be performed while fixing $d$.
\end{enumerate}
\end{lemma}

\begin{pf}
See Fig. \pr{f1}. 
\end{pf}

\begin{figure}[h]
\scalebox{0.6}{\includegraphics{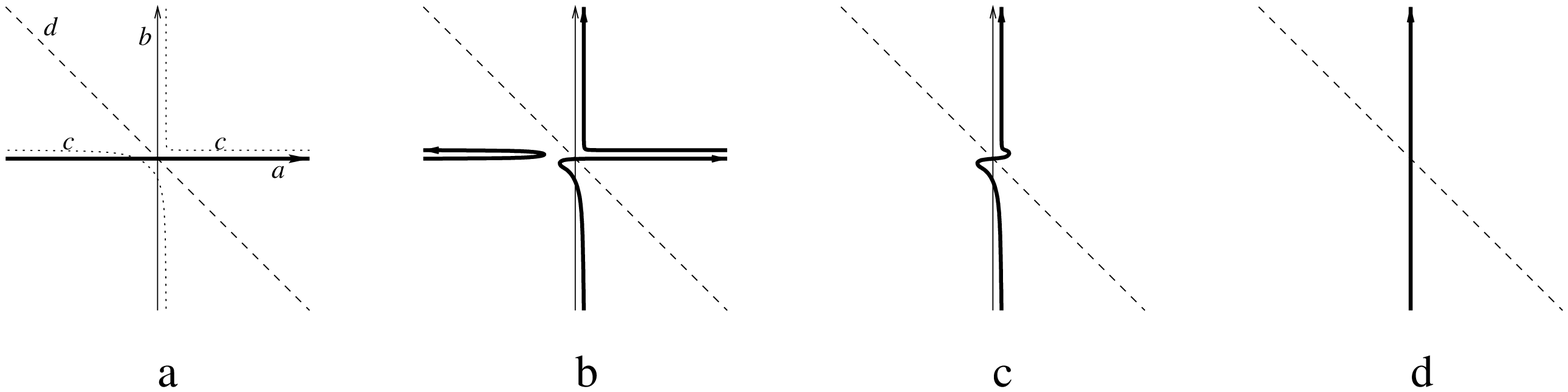}}
\caption{}\label{f1}
\end{figure}

\begin{figure}[h]
\scalebox{0.6}{\includegraphics{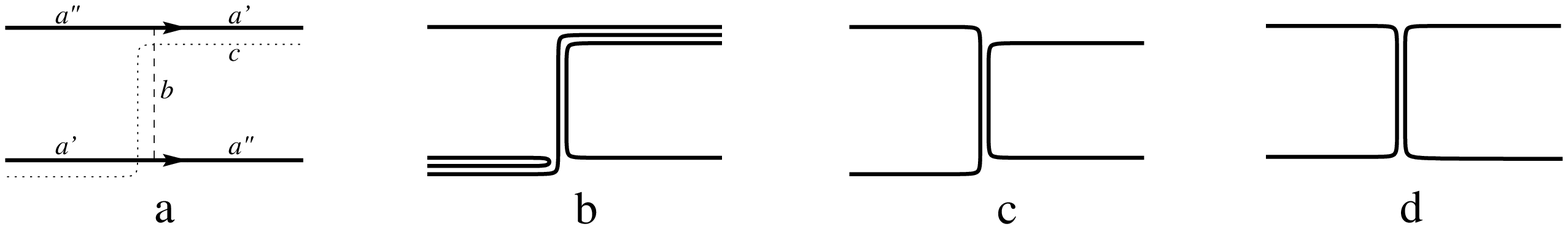}}
\caption{}\label{f2}
\end{figure}

\begin{lemma}\label{l2}
Let $a$ be a circle in $F$ and $b$ an arc connecting
two points of $a$, and whose interior is disjoint from $a$.
Assume that at the two endpoints of $b$, $a$ passes $b$
in the same direction (as in Fig. \pr{f2}a.)
Let
$a',a''$ be the two parts of $a$ into which it is separated by $\pa b$ 
and let $c=b\cup a'$. 
Then $(\T)^2$ followed by an isotopy performed in a thin 
neighborhood of $c$, maps $a$
onto the circle obtained by surgering $a$ along the arc
$b$ as in Fig. \pr{f2}d.
\end{lemma}

\begin{pf}
See Fig. \pr{f2}. 
\end{pf}

The setting for the following lemma is that of 
sections \pr{A},\pr{B}:

\begin{lemma}\label{l0}
Assume $(V,g)$ is not of the two special cases
$\dim V = 2$, $\A=0$ and $\dim V = 4$, $\A=0$.
Let $w_1, \dots ,w_k\in V$, $k\geq 0$, be independent vectors with
$g(w_i)=1$ and $B(w_i,w_j)=0$ for all $1\leq i,j \leq k$. Let 
$W=\langle w_1,\dots,w_k\rangle$ (the subspace spanned by $w_1,\dots,w_k$)
and let $a_1,a_2\in\wb - W$ be
two vectors with $g(a_1)=g(a_2)=1$ and $B(a_1,a_2)=0$. Then there exists
$c\in\wb$ with $g(c)=1$ and $B(a_1,c)=B(a_2,c)=1$. 
\end{lemma}

\begin{pf}
Assume first $k>0$.
We first find $b\in\wb$ such that $B(a_1,b)=B(a_2,b)=1$.
Since $a_1,a_2\not\in W=(\wb)^\bot$ there are
$b_1,b_2\in\wb$ with $B(a_1,b_1)=B(a_2,b_2)=1$. If also
$B(a_1,b_2)=1$ or $B(a_2,b_1)=1$ then we have a $b$. Otherwise $b_1 + b_2$
is our $b$. If $g(b)=1$ we are done with $c=b$, otherwise take $c=b+w_1$. 

Now assume $k=0$, so $W=\{ 0 \}$ and $\wb=V$.
If $a_1=a_2=a$, take $b\in V$ with $B(a,b)=1$. If $g(b)=1$ we are done with
$c=b$, otherwise define $U=\langle a,b \rangle$.
If $a_1\neq a_2$, take $b_1,b_2\in V$ with $B(b_1,b_2)=0$ and 
$B(a_i,b_j)=\delta_{ij}$ (Proposition \pr{pc3}). 
If $g(b_1+b_2)=1$ we are done with $c=b_1+b_2$, 
otherwise define $U=\langle a_1,a_2,b_1,b_2\rangle$.
In either case $B(x,y)$ is non-degenerate on $U$, and so 
it is non-degenerate on $U^\bot$.
If $\dim V > \dim U$ it follows that $g$ 
cannot be identically 0 on $U^\bot$. 
Take any element 
$d\in U^\bot$ with $g(d)=1$ then we are done with 
$c=b+d$ or $c=b_1+b_2+d$ respectively. 
So we are left with the case $\dim V = \dim U$.
If $\dim V = 2$ then we have assumed $\A=1$ and so we must have had $g(b)=1$.
If $\dim V = 4$ then again we have assumed $\A=1$ and so 
$g(b_1)\neq g(b_2)$ (since $g(a_1)=g(a_2)=1$)
so again we must have had $g(b_1+b_2)=1$.

\end{pf}

\begin{remark}\label{r3}
When $\dim V = 4$ then in the proof of Lemma \pr{l0} above, we haven't used 
the additional assumption that $\A=1$ in the following two cases:
(1) When $k>0$. (2) When $k=0$ and $a_1=a_2$ (since then $\dim V > \dim U$.)
\end{remark}

\begin{lemma}\label{l3}
Assume $F,g$ are not of the two special cases
$\g=1,\A=0$ and $\g=2,\A=0$. 
Let $w_1,\dots,w_k$ be disjoint circles in $F$ with $g([w_i])=1$
and such that $[w_1],\dots,[w_k]$ are independent in $\hh$
(which is equivalent to $\w$ not separating $F$.) 
Let $a_1,a_2$ be oriented circles in $F$ with $g([a_1])=g([a_2])=1$
and such that $a_1,a_2$ are each disjoint from $\w$ and
$[a_1],[a_2]\not\in\langle[w_1],\dots,[w_k]\rangle$.
Then there is a sequence $h_1,\dots,h_m$ of good maps of type
1 and 2 which all fix $\w$ and such that 
$h_1\circ\cdots\circ h_m$ (followed by an isotopy fixing $\w$) 
maps $a_1$ orientation preservingly onto $a_2$. 
\end{lemma}

\begin{pf}
Assume first that $[a_1]\cdot[a_2]=1$. If actually $|a_1\cap a_2|=1$ 
then we are done by Lemma \pr{l1} since the merge $c$
of $a_1$ and $a_2$ satisfies $g([c])=g([a_1])+g([a_2])+[a_1]\cdot[a_2]=1$
and so $\T$ is a good map of type 2.
So assume $|a_1\cap a_2|$ 
is some odd number $>1$. Then necessarily there exist two
consecutive crossings along $a_2$, at which $a_1$ crosses $a_2$ 
in the same direction. Applying the map $(\T)^2$ of Lemma \pr{l2} 
($a$ is here $a_1$ and $b$ is a portion of $a_2$) reduces 
$|a_1\cap a_2|$ by precisely 2, and so again 
$[a_1]\cdot[a_2]=1$ and so we may continue by induction.

Assume now $[a_1]\cdot[a_2]=0$. By Lemma \pr{l0}
there is $x\in\hh$ with $g(x)=1$, $x\cdot[a_1]=x\cdot[a_2]=1$
and $x\cdot [w_i]=0$ for all $i$. 
($x\not\in \langle [w_1],\dots,[w_k]\rangle$ follows.)
There exists a circle $c$
in $F$ with $[c]=x$ and such that $c$ is disjoint from each $w_i$.
(Start with any embedded representative and surger it along the $w_i$s
until it is disjoint from all of them. 
This is possible since the number of intersection points with each $w_i$
is even. Then connect the various components
to each other by surgery. This is possible since $F-\w$ is connected
and since there are no orientations to consider.) 
By the previous case we may now map $a_1$ onto $c$, and from there, 
orientation preservingly onto $a_2$. 
\end{pf}

\begin{remark}\label{r4}
If $F,g$ are of the special case $\g=2,\A=0$ then if 
in Lemma \pr{l3} we further assume 
either that $[a_1]\cdot[a_2]=1$ or that $[a_1]=[a_2]$ or that
$k>0$ then it follows from the proof of Lemma \pr{l3} and
from Remark \pr{r3}, that the conclusion of Lemma
\pr{l3} still holds.
\end{remark}

\begin{dfn}\label{um}
Let $\g=2,\A=0$. A map $\U\in\m$ such that $\U_*$ is a $U$-map
on $\hh$ (Definition \pr{u}) will again be called a $U$-map.
(Such maps clearly exist.)
\end{dfn}

\begin{thm}\label{t3}
If $F,g$ are not of the special case $\g=2,\A=0$ then 
$\m$ is generated by the good maps. 
In the mentioned special case, $\m$ is generated 
by the good maps and any one $U$-map.
\end{thm}

\begin{pf}
We first assume we are not in the two special cases appearing
in Lemma \pr{l3}, in particular, we are not in the special case
of this theorem. Let $n=\g$ and let $a_1,\dots,a_n,b_1,\dots,b_n$
be circles such that $|a_i\cap a_j|=|b_i\cap b_j|=0$ ($i\neq j$),
$|a_i\cap b_j|=\delta_{ij}$ and $g([a_i])=1$. 
(Start with such circles without the assumption on $g$. Then if 
for some $i$ both circles have $g=0$, replace one of them with their
merge. Then exchange $a_i$ with $b_i$ if necessary.)

Now let $h\in\m$. We will compose $h$ with good maps until (after isotopy) 
we arrive at the identity. We may first use Lemma \pr{l3}
to bring the $a_i$s one by one back to place. Indeed, after 
$a_1,\dots,a_k$ are in place, consider $a_1,\dots,a_k$, 
$h(a_{k+1})$, $a_{k+1}$, as the 
$w_1,\dots,w_k$, $a_1$, $a_2$ of Lemma \pr{l3}, respectively,
assigning orientations to $a_{k+1}$ and $h(a_{k+1})$ which
correspond via $h$.

So assuming $h$ fixes all $a_i$s, we bring each $b_i$
back to place, and assume we have already done this for all $i<j$.
Denote $a=a_j$, $b=b_j$ and let
$P$ be the intersection point of $a$ and $b$. Since $a$ is fixed
by $h$, $h(b)$ must also pass $P$, and since $h$ is orientation preserving,
$b$ and $h(b)$ must cross $a$ at $P$ in the same direction, 
(where orientations on $b$ and $h(b)$ correspond via $h$.)
Our permanent assumption that any circles we consider intersect transversally,
may still be maintained for the intersection of $b$ and $h(b)$ at $P$.
Assume first that $P$ is the only intersection point between $b$ and $h(b)$.
Let $c$ be the negative merge of $b$ and $h(b)$. $g([h(b)])=g([b])$ 
(since $h$ preserves $g$) and so $g([c])=g([b])+g([h(b)])+[b]\cdot[h(b)]=1$.
By Lemma \pr{l1}, $\T$ and an isotopy bring $h(b)$ orientation 
preservingly onto $b$, and the Dehn twist and isotopy may be performed 
while fixing $a$ and while fixing all other $a_i$s and all $b_i$s with
$i<j$. (Our $h(b)$, $b$, $c$ and $a$, correspond to $a$, $b$, $c$ and $d$
of Lemma \pr{l1} respectively.)

So assume now that there are additional 
intersection points between $h(b)$ and $b$ besides $P$. 
Choose a side of $a$ in $F$.
Let $X$ be the first additional intersection point,
when moving from $P$ along $b$ into the chosen side. 
Let $a'$ be a circle parallel
and close to $a$ on the chosen side. $g([a'])=1$.
If $a'$ is close enough to $a$, then $a'$ intersects both $b$ and $h(b)$ 
each at
a single point, and so $\TT_{a'}$ applied to $h(b)$ will add precisely
one intersection point between $h(b)$ and $b$. We denote this new 
intersection point by $P'$. If we choose the orientation of the Dehn twist 
$\TT_{a'}$ correctly, then $h(b)$ will cross $b$ at $P'$ in the same
direction that it crosses $b$ at $X$. 
Let $c$ be the circle which is the union of the subarcs of $b$ and $h(b)$ 
which are defined by $P'$ and $X$, and which do not contain $P$, and so
$c$ is disjoint from $a$. By Lemma \pr{l2}, $(\T)^2$ applied 
to $h(b)$ reduces the intersection between $h(b)$ and $b$ by precisely 2.
And so we have first increased the intersection by 1 and then decreased 
it by 2, and so we may continue by induction.

We are now 
in the situation that all $a_i$s and $b_i$s are fixed, so we are 
left with performing a map fixing all $a_i$s and $b_i$s, which is
equivalent to performing a boundary fixing map on a disc with $n$ holes. 
The group of all such maps is known to be
generated by Dehn twists. Now, any circle in the complement
of the $a_i$s and $b_i$s is bounding in $F$, and so these Dehn twists
are good maps of type 3. This completes the proof for the general case.

The case $\g=1,\A=0$ does not rely on the above, and will not be 
used in the sequel. We differ it to the end of Section \pr{D1}.

We are left with the case $\g=2,\A=0$.
We will show that by adding any $U$-map $\U$, 
the above proof can be made to work.
By Remark \pr{r4}, the only problem we have is 
when moving the first circle $a_1$. Let $V_1,V_2$ be the pair of subsets
of $\hh$ from the definition of $U$-map, and say $[a_1]\in V_1$.
If $[h(a_1)]\in V_2$ then $[\U\circ h(a_1)]\in V_1$,
and so we may assume $[a_1]$ and $[h(a_1)]$ are both in $V_1$. But then either 
$[a_1]=[h(a_1)]$ or $[a_1]\cdot [h(a_1)]=1$. 
By Remark \pr{r4} again, Lemma \pr{l3} applies, and so the 
above process works.

\end{pf}

\begin{remark}\label{r5}
If $p:\M\to GL(\hh)$ denotes the map $h\mapsto h_*$ then $p(\M)$ 
is known to be $Sp(\hh)$, the group of maps
preserving the intersection form on $\hh$.
Since $\oh\su Sp(\hh)$ and $\m = p^{-1}(\oh)$, $p|_{\m}:\m\to\oh$
is onto. Since $B(x,y)$ is unique up to isomorphism
for every given dimension $2n$, we see that our Theorem \pr{t3}
implies Theorem \pr{tc} (the Cartan-Dieudonne Theorem for the field $\C$.)
\end{remark}

\section{Quadruple Points of Regular Homotopies}\label{D}

Let $A$ be an annulus. There are two regular homotopy classes
of immersions $i:A\to\E$. This follows from the Smale-Hirsch Theorem ([H]) 
and the fact that $\pi_1(SO_3)=\C$. The first class is represented 
by an embedding whose image lies in a plane in $\E$, and the second class,
by an embedding  that differs from the former by one full twist.
For an immersion $i:A\to\E$ we let $G(i)\in\C$ be 0 or 1 according
to whether $i$ belongs to the first or second class,
respectively. 

\begin{dfn}
If $X\su Y$ then $N(X,Y)$ will denote a regular neighborhood 
of $X$ in $Y$.
\end{dfn}

Let $F$ be a closed orientable surface, and $i:F\to\E$
an immersion, then $i$ determines a quadratic form $g^i:\hh\to\C$ 
whose associated bilinear form is the intersection form on $\hh$, as follows:
Let $x\in\hh$, choose a circle $c$ in $F$ such that $[c]=x$ and 
define $g^i(x)=G(i|_{N(c,F)})$.
One needs to verify that $g^i(x)$ is independent of the choice of
$c$ and that $g^i(x+y)=g^i(x)+g^i(y)+x\cdot y$.
This has been done in [P] in the more general setting of
surfaces which are not necessarily orientable. 
(It is then necessary for the quadratic form to attain values
in $({1\over 2}\Z)/2=\Z/4$ rather than $\C$, 
to accommodate for the half twists 
of Mobious bands, and the Arf invariant attains values in $\Z/8$.)
For immersions $i,i':F\to\E$, $i\sim i'$ 
will denote that $i$ and $i'$ are regularly homotopic
in $\E$. The following has been shown in [P]:

\begin{thm}\label{Pi1}
Let $i,i':F\to\E$ be two immersions.
\begin{enumerate}
\item $g^i=g^{i'}$ iff $i\sim i'$. 
\item $\Ai =\Ar(g^{i'})$
iff there exists a diffeomorphism $h:F\to F$ such that $i\sim i'\circ h$.
\item $\Ar(g^e)=0$ for any embedding $e:F\to\E$.
\end{enumerate}
\end{thm}

Let $\hM$ be the group of \emph{all} 
diffeomorphisms of $F$ (not necessarily orientation preserving) 
up to isotopy.
Given a quadratic form $g$ on $\hh$ 
(whose associated bilinear form is the intersection form)
let $\hm$ be the subgroup of $\hM$ defined by
$\hm=\{h\in\hM: h_*\in\oh \}$.
($\m$ is then a subgroup of index 2 in $\hm$.)
Now let $i:F\to\E$ be an immersion and $h:F\to F$ a diffeomorphism. 
By Theorem \pr{Pi1}(1), $i\sim i\circ h$ 
iff $g^i=g^{i\circ h}$. It is easy to see that 
$g^{i\circ h}=g^i\circ h_*$ and so we get:

\begin{prop}\label{aa}
$i\sim i\circ h$ iff $h\in \hmi$.
\end{prop}

Let $H_t:F\to\E$ be a generic regular homotopy. We denote by $q(H_t)\in\C$
the number mod 2 of quadruple points occurring in $H_t$.
The following has been shown in [N]:

\begin{thm}\label{Q}
Let $F$ be any closed surface (not necessarily orientable or connected.)
If $H_t,G_t:F\to\E$ are two generic regular homotopies 
between the same two generic immersions, then $q(H_t)=q(G_t)$.
\end{thm}

\begin{dfn}\label{qii}
Let $i,i':F\to\E$ be two regularly homotopic generic immersions.
We define $Q(i,i')\in\C$ by $Q(i,i')=q(H_t)$, where $H_t$ is any 
generic regular homotopy between $i$ and $i'$. 
This is well defined by Theorem \pr{Q}. 
\end{dfn}

If $H_t, G_t:F\to\E$ are two regular homotopies
such that the final immersion of $H_t$ is the  
initial immersion of $G_t$, then $H_t * G_t$ will denote the regular 
homotopy that performs $H_t$ and then $G_t$.

\begin{lemma}\label{qh}
Let $i:F\to\E$ be a generic immersion. The map $\hmi\to\C$
given by $h\mapsto Q(i,i\circ h)$ is a homomorphism.
\end{lemma}

\begin{pf}
Let $h_1,h_2\in\hmi$ and let $H^k_t$ be a generic regular homotopy 
from $i$ to $i\circ h_k$, $k=1,2$. Then $H^1_t * (H^2_t\circ h_1)$
is a regular homotopy from $i$ to $i\circ h_2\circ h_1$ and
$q(H^1_t * (H^2_t\circ h_1))=q(H^1_t)+q(H^2_t)$.
\end{pf}

Recall that for $T\in \ov$ we have defined $\p(T)=\r(T-Id)\bmod{2}$, 
and have shown that $\p:\ov\to\C$ is a homomorphism
(Theorem \pr{t1}). For $h\in\hm$ let $\ep(h)\in\C$ 
be 0 or 1 according to whether $h$ is orientation preserving or reversing, 
respectively, and let $n=\g$.
Since $\ep:\hm\to\C$ and $h\mapsto h_*$ are also homomorphisms, 
the following $\hp:\hm\to\C$ is a homomorphism:
$$\hp(h) \ = \ \p(h_*)+(n+1)\ep(h) \ 
= \ \bigg(\r(h_*-Id)+(n+1)\ep(h)\bigg)\bmod{2}.$$

Our purpose is to show:

\begin{thm}\label{main1}
Let $i:F\to\E$ be a generic immersion. Then for any $h\in\hmi$:
$$Q(i,i\circ h)= \hp(h).$$ 
\end{thm}

Let $i:F\to\E$ be an immersion and let $c$ be a circle in $F$  
such that $c$ is disjoint from the multiplicity set of $i$.
Adding a \emph{ring} to $i$ along $c$ will mean to change $i$ into a 
new immersion $i'$ in the following way:
Let $U = N(i(c),\E)$, thin enough so that
$A=i^{-1}(U)$ is an annulus which is still disjoint from the multiplicity set.
Let $D$ and $a$ be a disc and an arc. Let $f_1:a\to D$ be a proper embedding
and let $f_2:a\to D$ be a proper immersion with one transverse intersection
and such that  $f_1|_{N(\pa a,a)}=f_2|_{N(\pa a, a)}$.
Parametrize $U$ and $A$ as $D\x S^1$ and $a\x S^1$ so that $i|_A:A\to U$
will be given by $f_1\x Id : a\x S^1 \to D\x S^1$. Now, the new immersion $i'$
will be given by $f_2\x Id$ on $A$, and will be the same as $i$ outside $A$.
There are basically two ways for adding a ring to $i$ along $c$, 
depending on what 
side of $A$ in $\E$ the ring will be facing
(which in turn depends on our choice of $f_2:a\to D$.)
If $i$ is an embedding, then $i(F)$ separates $\E$ into two pieces,
one compact and one non-compact. They will be denoted $\cc_i$ and $\nc_i$ 
respectively. And so if $i$ is an embedding then we have a natural 
way for distinguishing the two possibilities for 
adding a ring along a given circle $c$, 
namely, the ring is facing either $\cc_i$ or $\nc_i$.

Note that $f_1\x Id , f_2\x Id : A \to U$ are homotopic relative $\pa A$
(but not regularly homotopic.) And so if $N=N(i(F),\E)$ 
with $N\supseteq  U$, then $i$ and $i'$ 
are homotopic in $N$ (but in general not regularly homotopic.)

\begin{figure}[t]
\scalebox{0.6}{\includegraphics{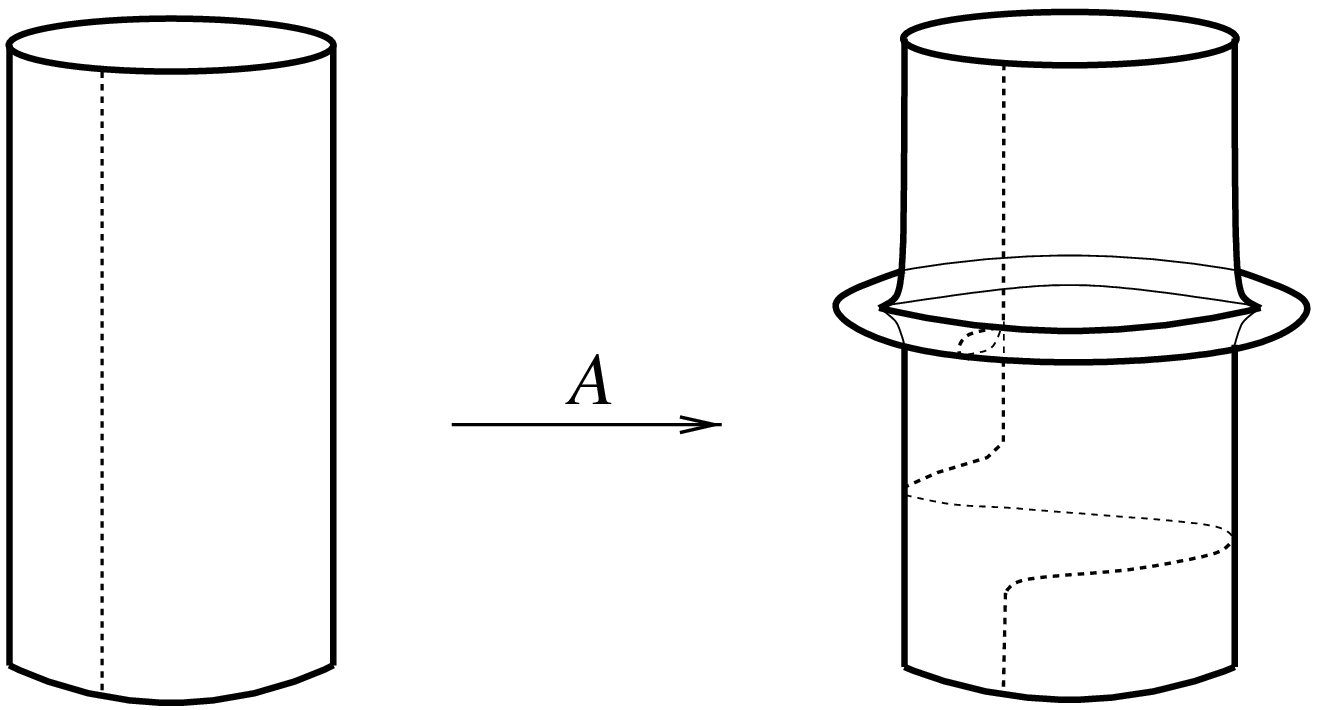}}
\caption{}\label{f3}
\end{figure}

\begin{figure}[t]
\scalebox{0.6}{\includegraphics{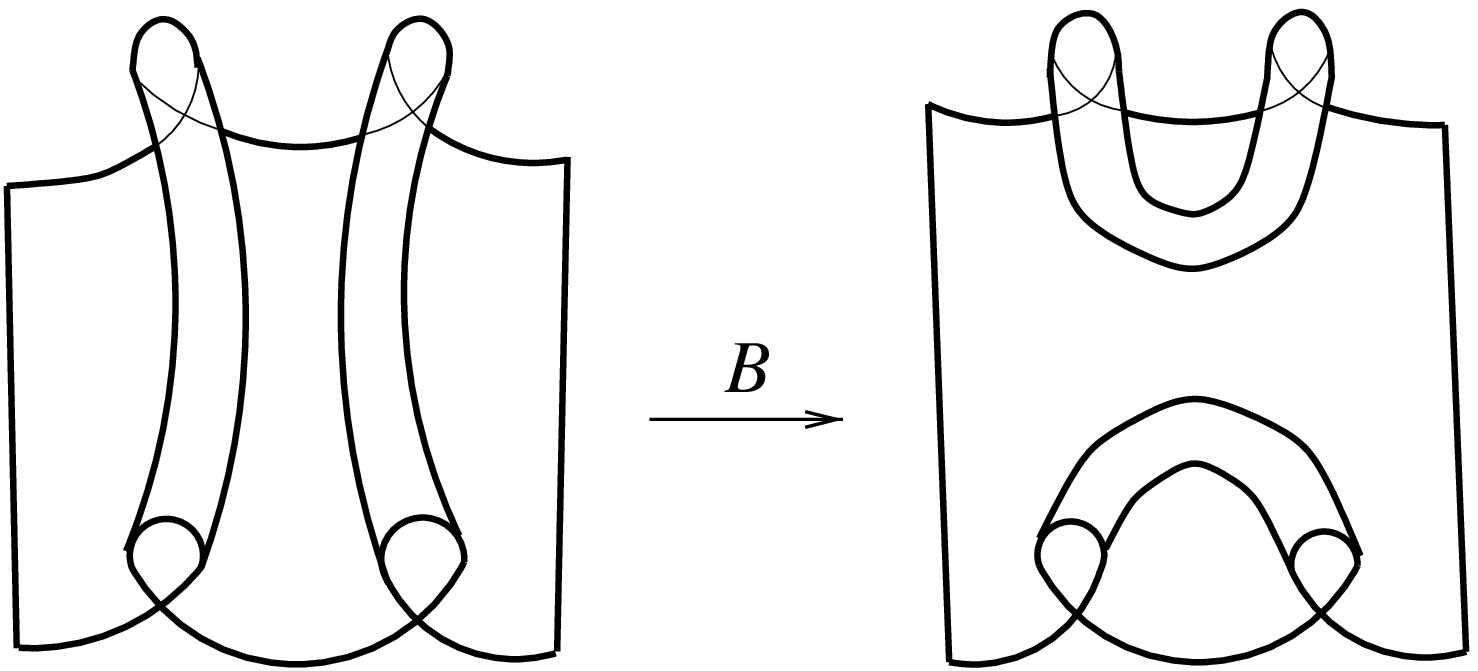}}
\caption{}\label{f4}
\end{figure}

We now present two moves on immersions,
that have been introduced in [N]. 
Let $S_0$ be an annulus and $S_1$ be a disc.
Move $A$ (resp. $B$) is a regular homotopy which is applied to a proper
immersion of $S_0$ (resp. $S_1$) into a ball $E$ and which fixes 
$N(\pa S_k , S_k), k=0,1$.
Move $A$ begins with the standard embedding of $S_0$ in $E$, 
and adds a ring and a Dehn twist along parallel essential circles in $S_0$.
The ring may face either side of $S_0$ in $E$ and the Dehn twist may
have either orientation. Fig. \pr{f3} shows one of the possibilities.
The reverse move, going from right to left in Fig. \pr{f3}, will be called
an $A^{-1}$ move.
Move $B$ is described in Fig. \pr{f4}. It begins with a specific 
immersion of $S_1$, with two arcs of intersection, and
replaces them with two other arcs of intersection. 
It is easy to see (and has been shown in [N])
that the initial and final immersions that we have presented
for the $A$ and $B$ moves, are indeed regularly
homotopic in $E$ (while keeping $N(\pa S_k, S_k)$ fixed.)
Move $A$ (resp. $B$) will be applied to an 
immersion $i:F\to\E$ inside a  
ball $E$ in $\E$ for which $i^{-1}(E)$ is an annulus (resp. disc)
and $i|_{i^{-1}(E)}:i^{-1}(E)\to E$ is as above.
(The rest of $F$ will be kept fixed.)
In particular, an $A$ move may be applied to a neighborhood of a circle
$c$ in $F$ iff $c$ is disjoint from the multiplicity set of $i$ and
there is an embedded disc $D$ in $\E$ such that $D\cap i(F)=i(c)$.
The move will then be performed in a thin $N(D,\E)$.
If the circle along which we perform an $A$ move happens to bound
a disc in $F$, then the Dehn twist that is produced is trivial,
and may be cancelled by rotating this disc.
The following has been shown in [N]:

\begin{prop}\label{pn1}
Let $S_0, S_1, E$ denote an annulus, disc and ball respectively.
\begin{enumerate}
\item For any generic
regular homotopy $A_t:S_0\to E$ that realizes an $A$ move, 
$q(A_t)=1$.
\item For any generic regular homotopy $B_t:S_1\to E$
that realizes the $B$ move, $q(B_t)=1$.
\end{enumerate}
\end{prop}

If $H_t:F\to\E$ is given by $H:F\x [0,1]\to\E$ then we denote by
$H_{-t}$ the regular homotopy given by $(x,t)\mapsto H(x,1-t)$.
Clearly $H_{-t}$ is generic iff $H_t$ is generic, and $q(H_t)=q(H_{-t})$.

\begin{lemma}\label{bb}
If $i\sim i'$ then $Q(i, i\circ h)=Q(i', i'\circ h)$
for any $h\in \hmi=\hM_{g^{i'}}$.
\end{lemma}

\begin{pf}
Let $J_t$ be a regular homotopy from $i$ to $i'$ and
$H_t$ a regular homotopy from $i$ to $i\circ h$.
Then $J_{-t}*H_t*(J_t\circ h)$ 
is a regular homotopy from $i'$ to $i'\circ h$
and $q(J_{-t}*H_t*(J_t\circ h))=q(H_t)$. 
\end{pf}

\begin{remark}\label{rbb}
After we prove Theorem \pr{main1}, we will know that
the assumption $i\sim i'$ in Lemma \pr{bb} is actually unnecessary,
as long as $h\in \hmi\cap \hM_{g^{i'}}$. This is so since $\hp(h)$ 
does not depend on $i$. (It is a function of $h$ only.)
\end{remark}

We begin the proof of Theorem \pr{main1}. First let $\g=0$.
By Lemma \pr{bb} (and since all immersions of $S^2$ in $\E$ are
regularly homotopic)
we may assume $i$ is an embedding onto the unit sphere. 
Let $r:\E\to\E$ be the reflection with respect to the $xy$-plane, and
$h\in\hM$ be such that $i\circ h =r \circ i$, then $h$ is
the unique non-trivial element of $\hM=\hmi$ and $\hp(h)=\ep(h)=1$. 
The following is a regular homotopy from $i$ to $i\circ h$: Perform an $A$
move on the equator of the sphere, such that the ring formed will be
facing $\nc_i$. Then exchange the northern and southern hemispheres,
arriving at $i\circ h$. By Proposition \pr{pn1} the $A$ move 
contributed 1 mod 2 quadruple points, then exchanging the hemispheres involved 
only double curves, and so $Q(i,i\circ h)=1$. 
We proceed by induction on $\g$ beginning with $\g=1$.
For the starting point $\g=1$, we will need to separate the cases 
$\Ai=0$ and $\Ai=1$.

\subsection{The case genus($F$)=1, Arf($g^i$)=0}\label{D1}

This has basically been done in [N]. We present it here with slight variation.
Let $T$ denote the torus. We will say an embedding $e:T\to\E$ is
\emph{standard} if its image is the torus $\tilde{T}\su\E$ 
obtained by rotating the 
circle $\{ y=0 , (x-2)^2 + z^2 =1 \}$ around the $z$-axis. 
Let $\tilde{m},\tilde{l}$ be the circles in $\tilde{T}$ given by
$\tilde{m}=\{ y=0 , (x-2)^2 + z^2 =1 \}$ and 
$\tilde{l}=\{ z=0 , x^2+y^2=1 \}$ and choose some orientations 
for $\tilde{m}$ and $\tilde{l}$.
For a standard embedding $e:T\to\E$, let $m_e$ and $l_e$ denote
the oriented circles in $T$ such that $e(m_e)=\tilde{m}$ and
$e(l_e)=\tilde{l}$.

Since $\Ai=0$ then
by Theorem \pr{Pi1}(2,3) $i$ is regularly homotopic to a standard embedding.
(Take an arbitrary standard embedding $i'$, then $i\sim i'\circ h$ for some 
diffeomorphism $h:T\to T$, but $i'\circ h$ is
again a standard embedding.) So by Lemma \pr{bb} we may assume $i$ itself 
is a standard embedding. 
By viewing $m=m_i,l=l_i$ as the basis for $H_1(T,\Z)$
(note $\Z$ coefficients) we identify $\hM$ with $GL_2(\Z)$. We will think
of any $h\in\hM$ both as a map from $F$ to $F$ and as a $2\x 2$ matrix. 
If $h\in\hM$ then $h_*:\hh\to\hh$ (now $\C$ coefficients)
presented with respect to the basis 
$[m],[l]$ is simply the $\C$ reduction of the matrix $h$.
$g^i([m])=g^i([l])=0$ and for $x=[m]+[l]$, $g^i(x) = 1$.
A matrix $h\in GL_2(\Z)$ is in $\hmi$ if $h_*$ preserves $g^i$.
This will happen iff $h_*(x)=x$ and $h_*(\{[m],[l]\})=\{[m],[l]\}$, 
i.e. iff the $\C$ reduction of $h$ is either 
$I=\left[\begin{smallmatrix} 1 & 0 \\ 0 & 1 \end{smallmatrix}\right]$
or $J=\left[\begin{smallmatrix} 0 & 1 \\ 1 & 0 \end{smallmatrix}\right]$.
By means of row and column operations one can show that this subgroup
of $GL_2(\Z)$ is generated by the following four elements:
$A_1=\left(\begin{smallmatrix} 1 & 2 \\ 0 & 1 \end{smallmatrix}\right)$,
$A_2=\left(\begin{smallmatrix} 1 & 0 \\ 2 & 1 \end{smallmatrix}\right)$,
$A_3=\left(\begin{smallmatrix} -1 & 0 \\ 0 & 1 \end{smallmatrix}\right)$,
$A_4=\left(\begin{smallmatrix} 0 & 1 \\ 1 & 0 \end{smallmatrix}\right)$.

Since the $\C$ reduction of $A_1,A_2,A_3$ is $I$ and that of $A_4$ is
$J$, and since $\p(I)=0$, $\p(J)=1$ and $n=1$,
we have $\hp(A_1)=\hp(A_2)=\hp(A_3)=0$ and $\hp(A_4)=1$.
We will now show that the values of $Q(i,i\circ A_k)$, $k=1,\dots,4$ 
are the same.

\begin{figure}[t]
\scalebox{0.6}{\includegraphics{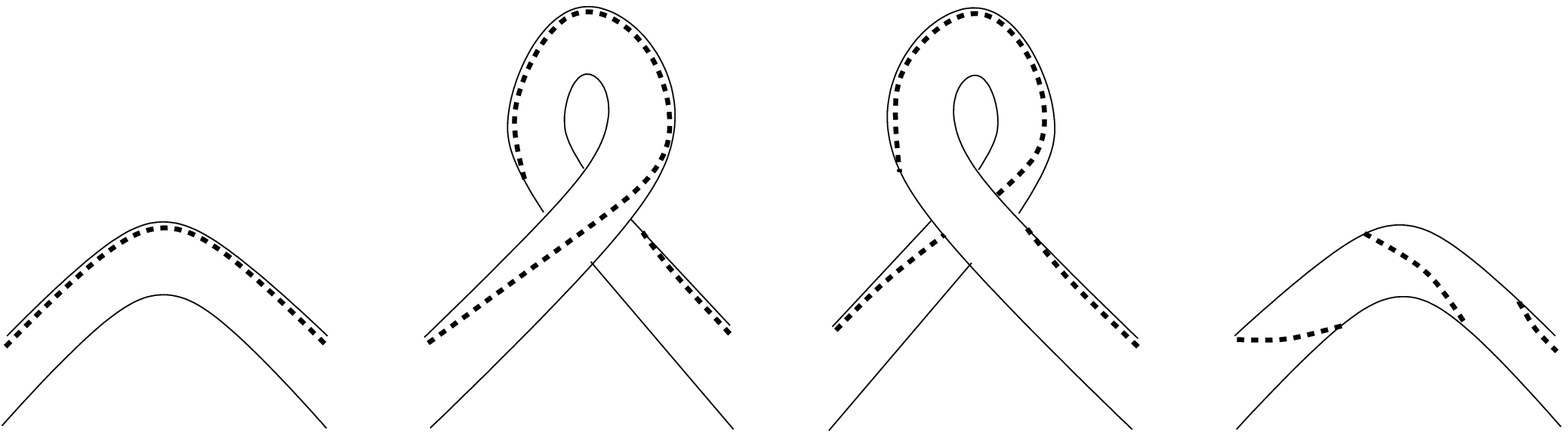}}
\caption{}\label{f5}
\end{figure}

$A_1$ and $A_2$ are $(\TT_m)^2$ and $(\TT_l)^2$ respectively.
Let $D$ be a compressing disc for $i(m)=\tilde{m}$ in $\E$ and $B=N(D,\E)$
thin enough so that $i(T)\cap B$ is a standard annulus in $B$. 
In $B$ we isotope $i(T)\cap B$ to be a thin tube, we then 
perform the ``belt trick'' on this tube 
(Fig. \pr{f5})
and then isotope $i(T)\cap B$ 
back to place. The effect of this regular homotopy is precisely $(\TT_m)^2$
and it involves only double curves, and so $Q(i,i\circ A_1)=0$. 
In the same way $Q(i,i\circ A_2)=0$.

Up to isotopy in $T$,
$i\circ A_3=r\circ i$ where $r$ is the reflection of $\E$ with respect to the 
$xy$-plane. This may be achieved by a regular homotopy similar
to the one we have used for the case $\g=0$, as follows:
Perform an $A$ move on each of 
the two circles $\{ z=0,x^2+y^2=1 \}$ and $\{ z=0,x^2+y^2=3 \}$, 
so that the ring 
formed by each of them is facing $\nc_i$ and such that the two Dehn twists
formed will have opposite orientations and so will cancel each other.
(The $A$ moves are performed inside thin neighborhoods of compressing discs
for the two circles.)
So we remain with just the two rings. 
We may now exchange the upper and lower halves of $T$ until we arrive at 
$i\circ A_3$. The two $A$ moves each
contributed $1\bmod{2}$ quadruple points and the final stage 
involved only double curves, and so all together $Q(i,i\circ A_3)=0$.

For $A_4$, isotope $T$ until it has the shape of a large sphere with 
a tiny handle at its north pole. Now exchange the northern and southern 
hemispheres. This will involve only double curves, and will result in
a sphere having a tiny handle at the south pole and a ring along 
the equator. We cancel this ring with an $A^{-1}$ move in a thin 
neighborhood of the plane of the equator, resulting in an embedding again. 
We may think of $m$ and $l$ as being contained in the tiny handle, and so
tiny compressing discs for $m$ and $l$ may be pulled along with 
the regular homotopy. The compressing disc of $m$ now lies 
in $\nc_i$, and the compressing disc
of $l$ now lies in $\cc_i$.
It follows that the final embedding may be isotoped to a standard
embedding i.e. to a map of the form $i\circ h$, and this $h$ is orientation
reversing and exchanges $m$ and $l$. And so $h$ must be either
$\left(\begin{smallmatrix} 0 & 1 \\ 1 & 0 \end{smallmatrix}\right)$ or
$\left(\begin{smallmatrix} 0 & -1 \\ -1 & 0 \end{smallmatrix}\right)$.
These two maps composed with $i$ are isotopic in $\E$ (via a half 
revolution about the $x$-axis) and so we may assume 
$h=\left(\begin{smallmatrix} 0 & 1 \\ 1 & 0 \end{smallmatrix}\right) =A_4$.
Our regular homotopy had a portion involving just double curves, then
one $A^{-1}$ move, and finally some isotopy, and so $Q(i,i\circ A_4)=1$.

This completes the proof that $Q(i,i\circ h)=\hp(h)$ for
every $h\in\hmi$ for $F$ a torus and $\Ai=0$. We now give the promised 
completion of the proof of Theorem \pr{t3}. We need to show that
$\m$ is generated by good maps when $\g=1,\A=0$.
Choose two oriented circles $a,b$ with $g([a])=g([b])=0$ 
and $|a\cap b|=1$ as a basis for
$H_1(T,\Z)$, thus identifying $\M$ with $SL_2 (\Z)$. Again we see
$h\in SL_2(\Z)$ is in $\m$ iff its $\C$ reduction is 
either $I$ or $J$.
By row and column operations, we then see that
$\m$ is generated by
$A_1=\left(\begin{smallmatrix} 1 & 2 \\ 0 & 1 \end{smallmatrix}\right)$,
$A_2=\left(\begin{smallmatrix} 1 & 0 \\ 2 & 1 \end{smallmatrix}\right)$ and
$A'=\left(\begin{smallmatrix} 0 & 1 \\ -1 & 0 \end{smallmatrix}\right)$.
Now $A_1$ and $A_2$ are $(\TT_a)^2$ and $(\TT_b)^2$. If
$c$ is the positive merge of $a$ and $b$ then $g([c])=1$
and $\T$ is
$\left(\begin{smallmatrix} 0 & 1 \\ -1 & 2 \end{smallmatrix}\right)$.
Since 
$\left(\begin{smallmatrix} 0 & 1 \\ -1 & 0 \end{smallmatrix}\right)=
\left(\begin{smallmatrix} 0 & 1 \\ -1 & 2 \end{smallmatrix}\right) 
\left(\begin{smallmatrix} 1 & 2 \\ 0 & 1 \end{smallmatrix}\right)$
we see $\m$ is generated by $(\TT_a)^2$, $(\TT_b)^2$ and $\T$.

\subsection{The case genus($F$)=1, Arf($g^i$)=1}\label{D2}

By Theorem \pr{Pi1}(1), $i$ is regularly homotopic to an immersion which 
is obtained from a standard embedding $e:T\to\E$ by adding a ring along 
the circle $c$ which is the positive merge of $m_e$ and $l_e$, 
and such that the ring is facing $\nc_e$. 
(One checks directly that such an immersion $i'$ has $\Ar(g^{i'})=1$,
but then $g^i=g^{i'}$ since on $V$ of dimension 2 there is only one
$g$ with $\A=1$.)
By Lemma \pr{bb} we may assume $i$ itself is this new immersion 
(Fig. \pr{f6}a).
Since all non-zero elements in $\hh$ have $g^i=1$,
it follows that $\ohi=GL(\hh)$ and so
$\hmi=\hM$. Since $\C$ is abelian, it is enough to verify
$Q(i,i\circ h)=\hp(h)$ only on normal generators, and we claim that
$B_1=\left(\begin{smallmatrix} -1 & 2 \\ 0 & 1 \end{smallmatrix}\right)$
and $B_2=\left(\begin{smallmatrix} 0 & 1 \\ 1 & 0 \end{smallmatrix}\right)$
are normal generators of $\hM$ where the identification with $GL_2(\Z)$
is via $m_e,l_e$. (There are no $m_i,l_i$ since such are defined only for
standard embeddings.) Indeed, 
$\left(\begin{smallmatrix} 0 & 1 \\ 1 & 0 \end{smallmatrix}\right)
\left(\begin{smallmatrix} -1 & 2 \\ 0 & 1 \end{smallmatrix}\right)
=\left(\begin{smallmatrix} 0 & 1 \\ -1 & 2 \end{smallmatrix}\right)$
which as we have noticed in the end of section
\pr{D1}, is a Dehn twist along $c$ (the positive merge of $m_e$ and $l_e$.) 
Now, any Dehn twist is a normal generator of $\M$, and since
$B_1, B_2$ are orientation reversing, they normally generate the whole of 
$\hM$. As above, we see that $\hp(B_1)=0$ and $\hp(B_2)=1$.

\begin{figure}[t]
\scalebox{0.6}{\includegraphics{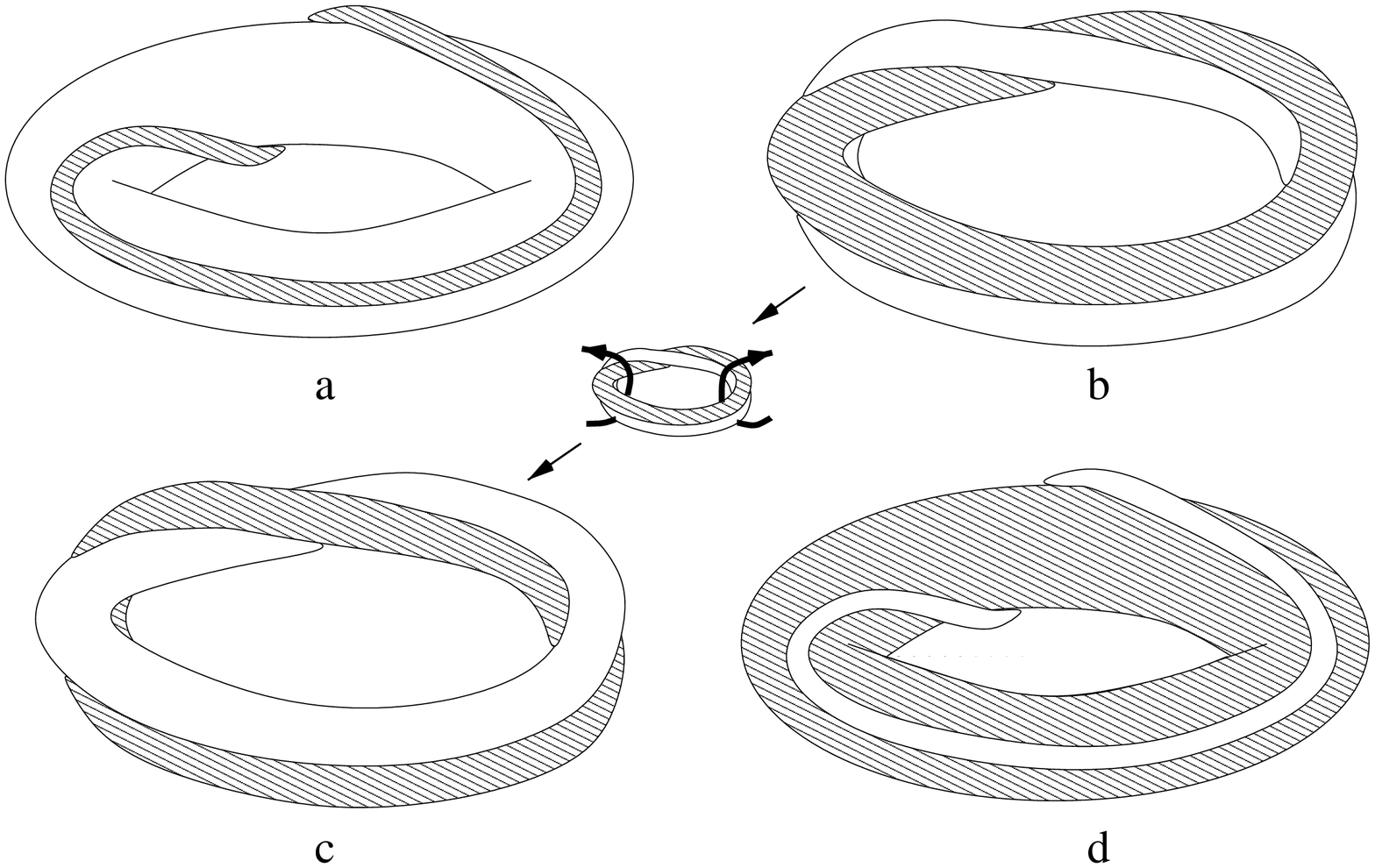}}
\caption{}\label{f6}
\end{figure}

The regular homotopy we construct for $B_1$ is as follows: Let the ring become 
thicker, and at the same time let the ``body'' of the torus become thinner, 
until they are equal in width, (Fig. \pr{f6}a$\to$b.) 
Now using the intersection circle as an axis of rotation,
we perform a half revolution of the torus around it,
interchanging the two equal-width rings, 
(Fig. \pr{f6}b$\to$c.)
We then return to our 
original position by reversing our first step. 
(Fig. \pr{f6}c$\to$d.) 
We thus arrive at an immersion of the form $i\circ h$. The map $h:T\to T$
may best be understood by looking at the intermediate stage of the regular
homotopy, Fig. \pr{f6}b$\to$c.
We see here that $h$ maps $m_e$ to $-m_e$ and $c$ to itself. 
i.e. the column $(1,0)^t$ to $(-1,0)^t$ and $(1,1)^t$ to itself.
And so indeed $h=B_1$. We had no quadruple points 
(actually no singular occurrences at all) and so $Q(i,i\circ B_1)=0$.

\begin{figure}[t]
\scalebox{0.6}{\includegraphics{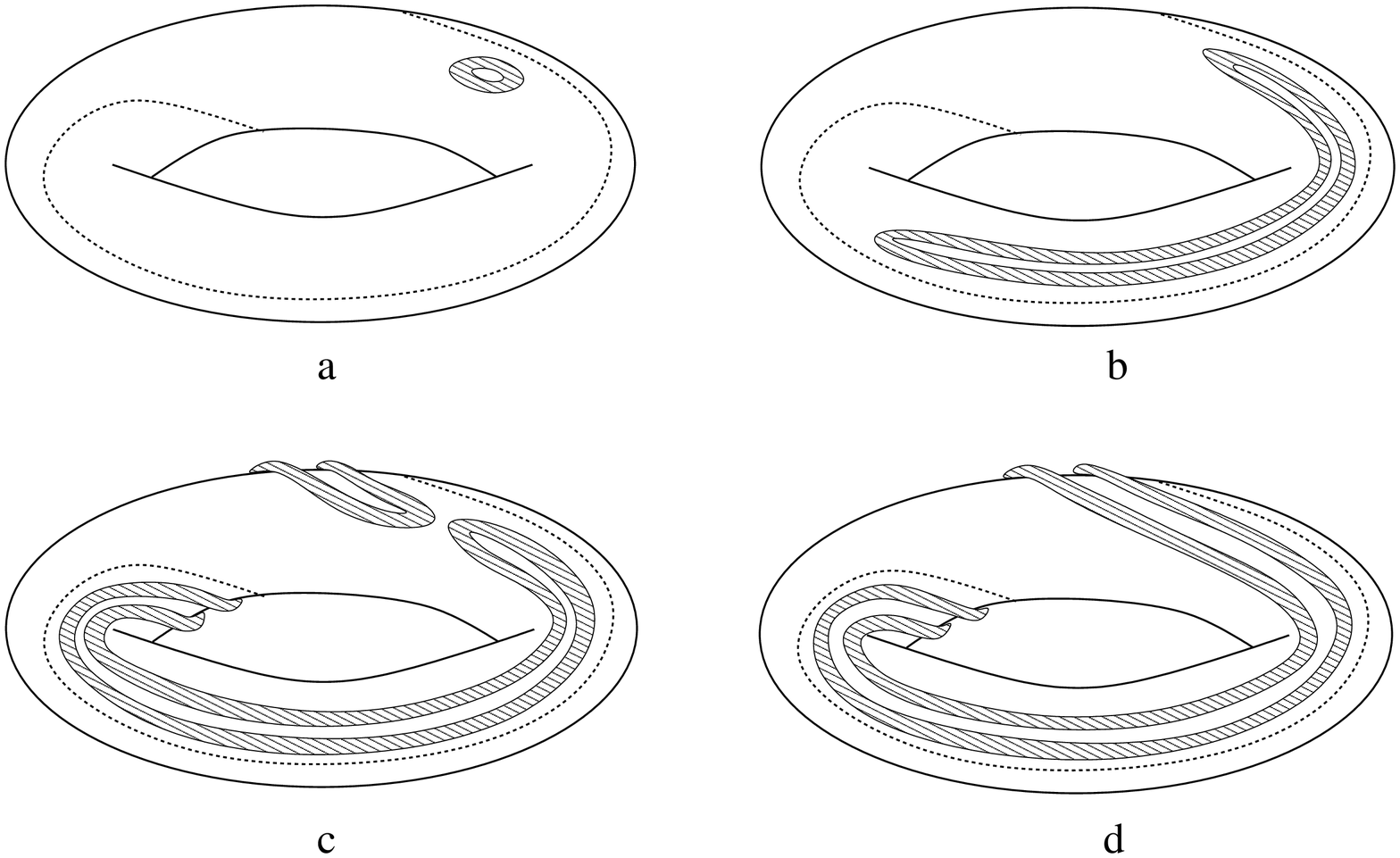}}
\caption{}\label{f7}
\end{figure}

For $B_2$, we first imitate the regular homotopy we have had for $A_4=B_2$ of  
Section \pr{D1} i.e. we perform that regular homotopy on $e$ and carry the
ring along. If before exchanging the upper and lower hemispheres
we make sure that the ring is situated at the tiny handle, 
then this exchange will have at most triple points, and 
the ring will not interfere with the $A^{-1}$ move, and so 
at the end of this process we will have $q=1$. 
The immersion $j$ we arrive at, is the immersion obtained by adding a ring $R$
to the embedding $e\circ B_2$ along the circle $e\circ B_2 (c) = e(c)$,
that is the same circle along which $R$ was originally situated,
but now $R$ is facing $\cc_e$ instead of $\nc_e$
and so $j$ is not of the form $i\circ h$. 
We fix this with a regular homotopy as follows:
Perform an $A$ move on a little disc bounding circle in $T$ near $R$,
forming a ring $R'$ facing $\nc_e$. See Fig. \pr{f7}a.
(The dotted line in Fig. \pr{f7} is the intersection curve
of $R$. $R$ itself is not seen since it is facing $\cc_e$.)
We then elongate $R'$ along side $R$, 
until it approaches itself from all the way around. 
(Fig. \pr{f7}a$\to$b$\to$c.)
We then perform a $B$ move, turning $R'$ into two rings which are parallel 
to $R$, but facing $\nc_e$. 
(Fig. \pr{f7}c$\to$d.) 
It is then easy to 
construct an explicit regular homotopy so that $R$ and 
the new ring which is adjacent to it, will cancel each other,
and with no quadruple points at all. 
(The idea is as follows: Let $f:a\to D$ be a proper immersion of an arc
$a$ into a disc $D$ with two loops facing opposite sides,
as in the front disc of Fig. \pr{f8}a.
There is a regular homotopy $f_t:a\to D$ fixing $N(\pa a,a)$,
from $f$ to an embedding as in the front disc of Fig. \pr{f8}b.
If $f_t$ is generic then it has at most triple points. 
Now the regular homotopy $f_t\x Id : a\x S^1 \to D\x S^1$ is a regular homotopy
which begins with a pair of rings facing opposite sides, and cancels them.
Fig. \pr{f8} depicts a portion of $a\x S^1$ in the initial and
final immersions.
Indeed since $f_t$ has at most triple points, so will $f_t\x Id$,
but $f_t\x Id$ is not generic. Nevertheless, any specific $f_t\x Id$
may serve as a guide for constructing an explicit generic regular homotopy 
with no quadruple points, and which begins and ends with the same immersions.)

\begin{figure}[t]
\scalebox{0.6}{\includegraphics{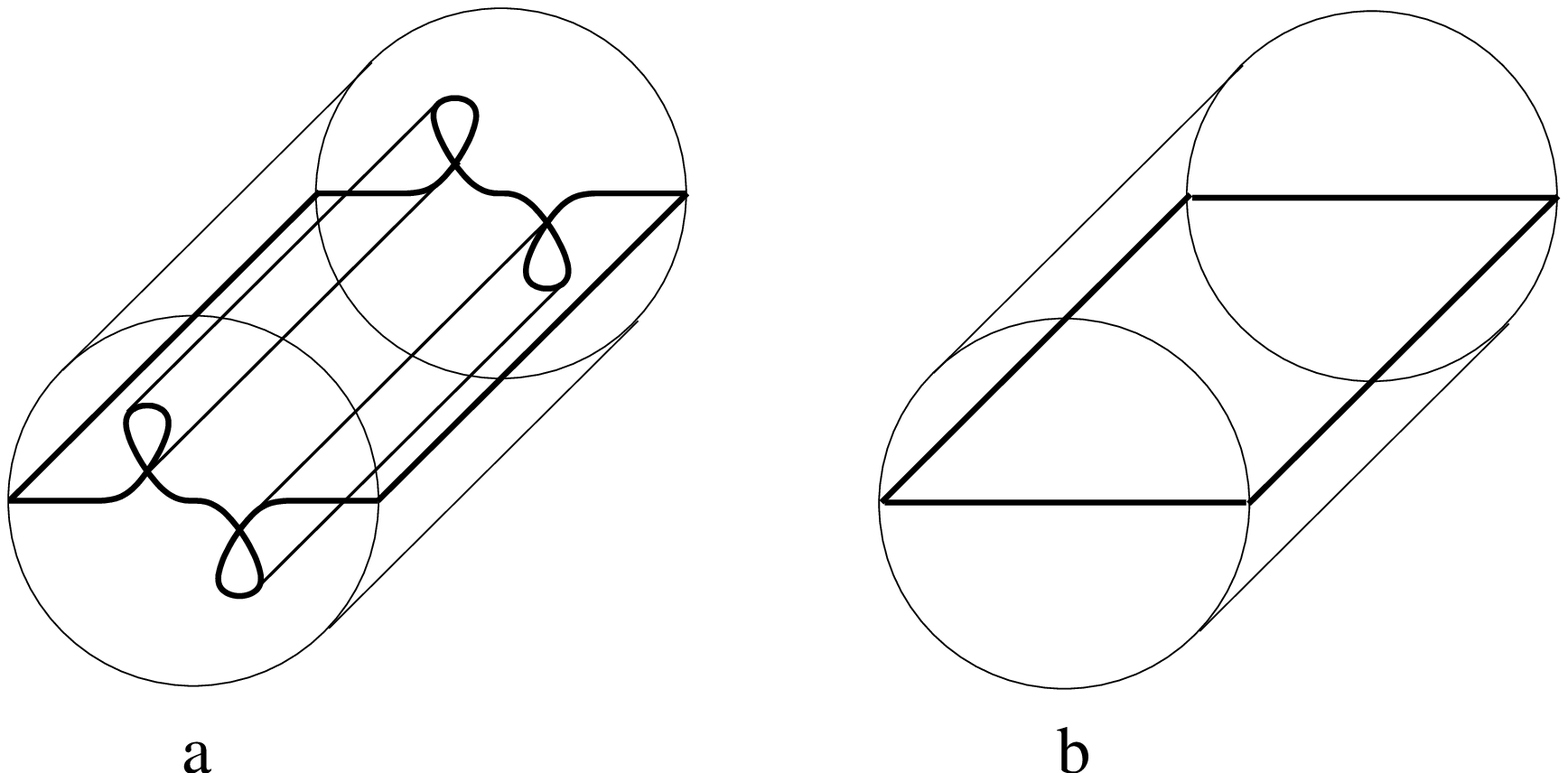}}
\caption{}\label{f8}
\end{figure}

Since $R$ and one of the new rings have disappeared, 
we are left with one ring facing $\nc_e$ and situated along
a circle parallel to $e(c)$. By pushing it precisely to $e(c)$ we finally
get an immersion of the form $i\circ h$.
The regular homotopy which started with $j$ and replaced 
the ring $R$ with a ring facing $\nc_e$,
took place inside some $N=N(e(F),\E)$, and so the homotopy class into
$N$ is still that of $e\circ B_2$ and so (up to isotopy in $F$)
the new immersion is indeed $i\circ B_2$. 
($i:F\to N$ is a homotopy equivalence, and so two 
diffeomorphisms $h,h':F\to F$ are isotopic in $F$ 
iff $i\circ h, i\circ h': F\to N$ are homotopic in $N$.)
Finally, our regular homotopy from $i$ to $j$ involved 
$1\bmod{2}$ quadruple points, then from $j$ to $i\circ B_2$ we had
one $A$ move, one $B$ move, 
and a regular homotopy with no quadruple points, and so all
together indeed $Q(i,i\circ B_2)=1$.

\subsection{The general case}\label{D3}
Assume $\g>1$.
By Theorem \pr{t3} $\mi$ is generated by Dehn twists and
squares of Dehn twists, and in the special case $\g=2,\Ai =0$
we also need a $U$-map. Other than the special $U$-map generator, which
will be dealt with last, each generator $h$
fixes all but a regular neighborhood of a circle $a$
(and we may assume $a$ does not bound a disc in $F$.) 
If $a$ is non-separating
then there is a circle $b$ in $F$ with $|a\cap b|=1$. $N(a\cup b,F)$
is a punctured torus, and so $c=\pa N$ is a circle separating
$F$ into two subsurfaces $F_1,F_2$ of smaller genus than $F$ and
with $h(F_k) = F_k$. If $a$ is separating,
then a nearby parallel circle $c$ will again separate $F$ in this way.
Now $\hmi$ needs one additional generator. Choosing any separating 
circle $c$ in $F$ there is clearly an orientation reversing $h:F\to F$
which preserves the two sides of $c$ and which induces the identity on $\hh$ 
and so $h\in\hmi$. And so finally we have a set of generators
for $\hmi$, each of which (except for the $U$-map)
preserves such a separation of $F$ by a circle $c$ into $F_1,F_2$ of
smaller genus.

\begin{figure}[t]
\scalebox{0.6}{\includegraphics{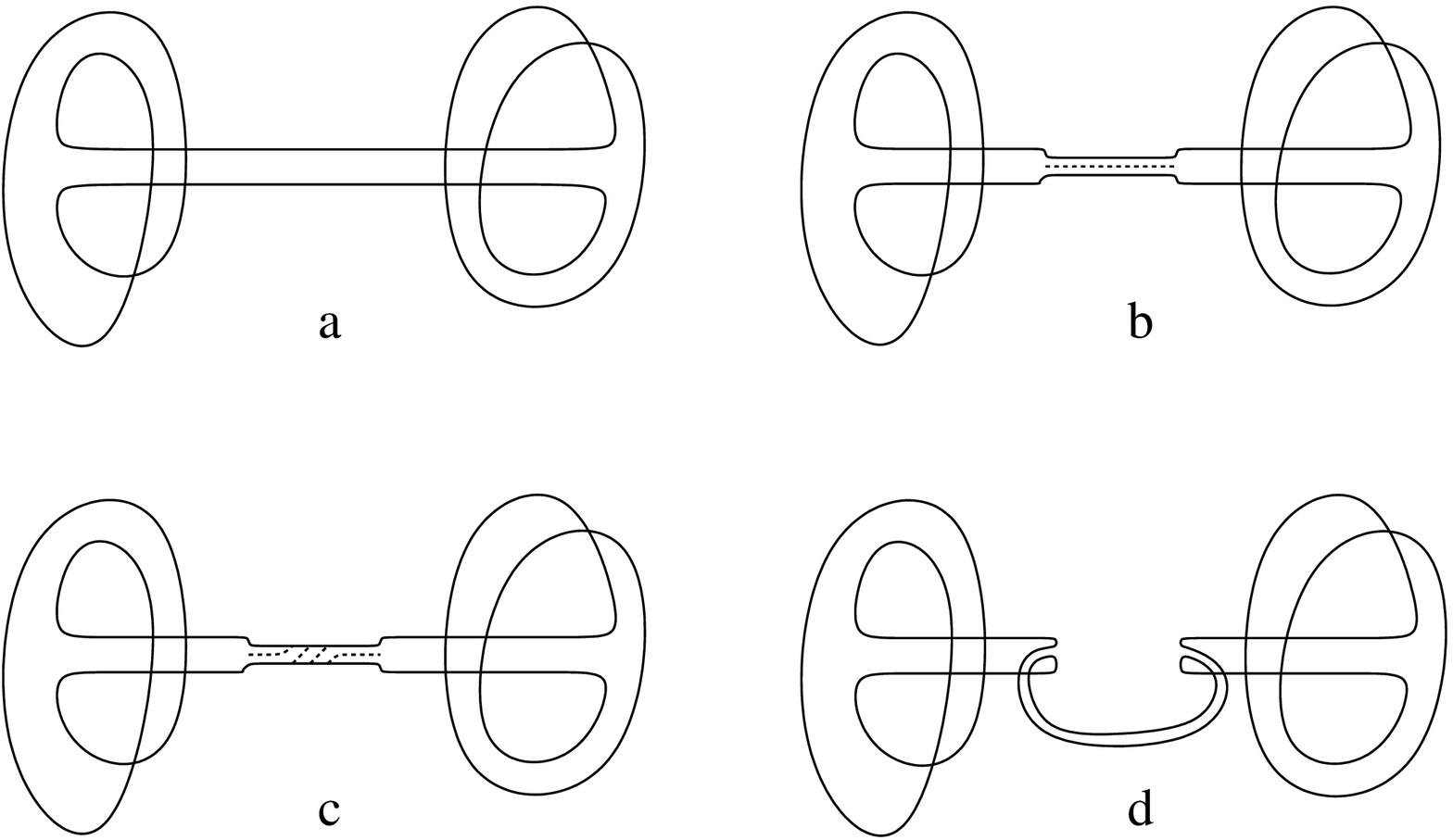}}
\caption{}\label{f9}
\end{figure}

Let $A=N(c,F)$. Slightly diminishing $F_1,F_2$ to be the
components of $F-intA$, we may still assume $h(F_k)=F_k$, $k=1,2$.
Since $c$ is separating in $F$, $g^i([c])=0$ and so $i|_A$ is regularly 
homotopic to a standard embedding of $A$, in the shape of a thin tube.
By means of [S] (namely the proof of Theorem 2.1) 
we may extend such a regular homotopy of $A$ to the whole of $F$. 
We now stretch this tube to be very long, 
at the same time pulling $F_1$ and $F_2$ rigidly away from each other until
they are disjoint. See Fig. \pr{f9}a. By taking a smaller
$A$ if necessary, we may assume $i(A)$ is disjoint from $i(F-A)$,
Fig. \pr{f9}b. By Lemma \pr{bb}, we may assume that this
is in fact our immersion $i$.
Let $\bar{F_1},\bar{F_2}$ be the closed surfaces
obtained by gluing a disc $D_k$ to $F_k$ and 
let $h_k:\bar{F_k}\to\bar{F_k}$ be an extension of $h|_{F_k}:F_k\to F_k$.
If the tube $i(A)$ is very thin,
then there is also a naturally defined extension 
$i_k:\bar{F_k}\to\E$ of $i|_{F_k}$.
We may further assume that the thin ball $B$ in $\E$ 
which is bounded by the sphere $i_1(D_1)\cup i(A) \cup i_2(D_2)$, 
is disjoint from $i(F-A)$. 

Since $h|_{F_k}$ preserves $g^i|_{H_1(F_k,\C)}$ then 
$h_k$ preserves $g^{i_k}$. It follows that there is a regular homotopy 
$H^k_t$ between $i_k$ and $i_k\circ h_k$. 
We perform $H^1_t$ and $H^2_t$ inside disjoint balls, and we let the
thin tube $A$ be carried along. If we make sure no quadruple points
occur in $D_1$ and $D_2$, and the thin tube $A$ does not pass triple points, 
then the regular homotopy $H_t$ induced on $F$ in this way will have the
sum of the number of quadruple points of $H^1_t$ and $H^2_t$. 

Now if $h$ is orientation preserving then so are $h_k$, in particular 
$h_k|_{D_k}$ is orientation preserving. So if we had carried the 
thin ball $B$ along with the tube $A$, then it would now approach the $D_k$s
from the same side it had for $i$. And so we may continue the regular 
homotopy on the tube $A$, still not passing through triple points, 
and cancelling all knotting by
having the thin tube pass itself, until it is back to its original place, 
and this will not contribute any quadruple points. 
However, the new embedding of $A$ may differ from $i\circ h |_A$  
by some number of Dehn twists as in Fig. \pr{f9}c.
We may resolve this by rigidly rotating say $F_1$ around the axis 
of the tube.

\begin{figure}[t]
\scalebox{0.6}{\includegraphics{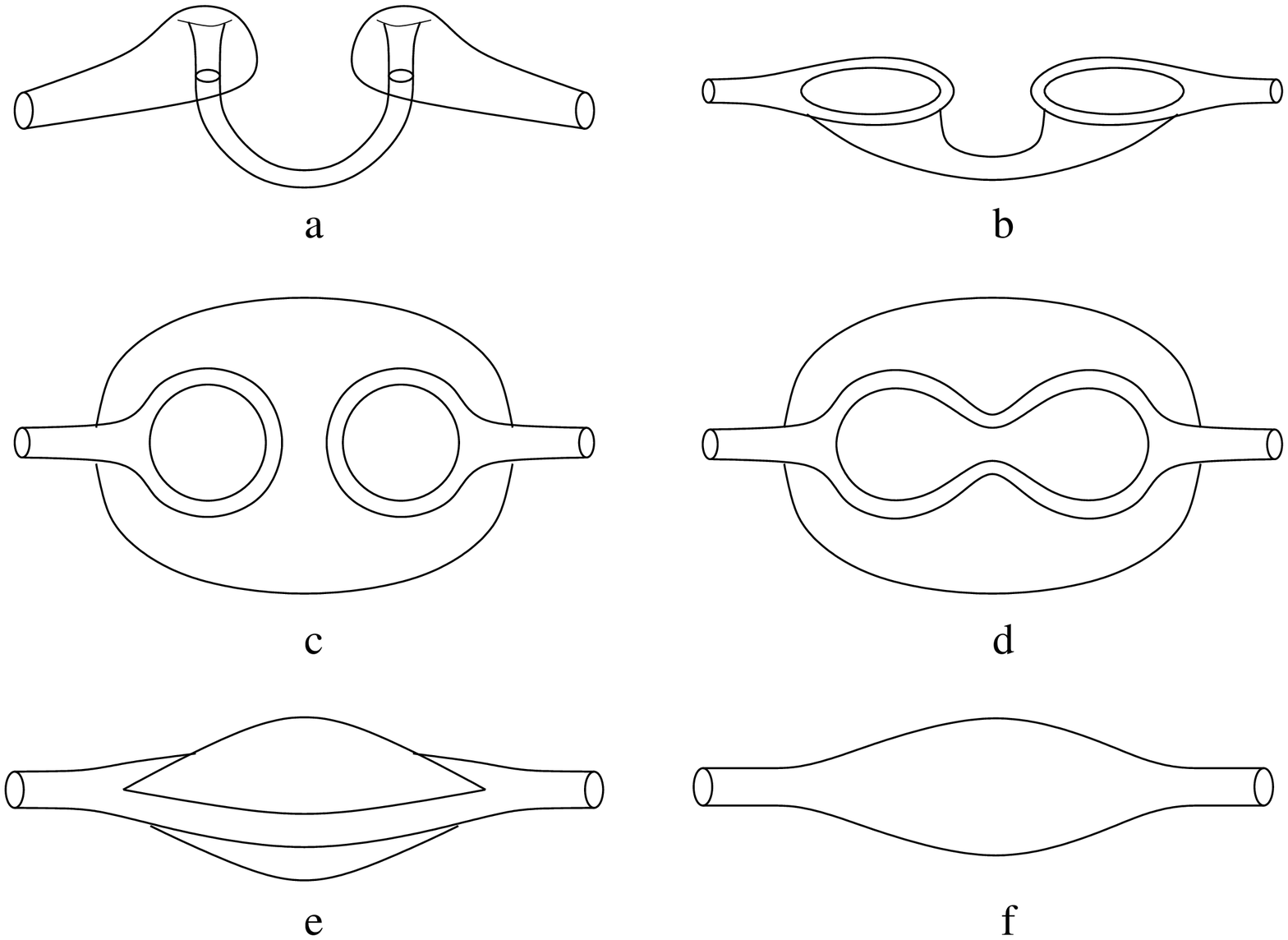}}
\caption{}\label{f10}
\end{figure}

If on the other hand $h$ is orientation reversing, then 
after applying $H^1_t$ and $H^2_t$ and carrying the tube along, 
the thin ball $B$ will approach both $D_k$s from the wrong side. 
And so after we cancel all knotting, the tube $A$ will be as in Fig \pr{f9}d.
Fig. \pr{f10} presents a regular homotopy that resolves this, and has 
$1\bmod{2}$ quadruple points. Fig. \pr{f10}a depicts the relevant part of 
Fig. \pr{f9}d, where the regular homotopy will take place. 
Fig. \pr{f10}a$\to$b$\to$c is a regular homotopy with no singular 
occurrences, or alternatively may be thought of as an ambient isotopy of $\E$.
It shows that we may view the immersion of $A$ as a sphere with two rings
facing outward, each of which has a tube attached to it. 
We now perform a $B$ move which joins the two rings into one ring with 
two tubes attached to it, Fig. \pr{f10}c$\to$d. 
Again by ambient isotopy, the ring may be brought to the equator, 
Fig. \pr{f10}d$\to$e. Finally we exchange the northern and 
southern halves of the sphere, arriving at an embedding, 
Fig. \pr{f10}e$\to$f. This regular homotopy involved a $B$ move and a
portion involving only double curves, and so indeed it had $1\bmod{2}$
quadruple points. We then continue to bring $A$ back to place. 
As above, the new embedding of $A$ may differ from $i\circ h|_A$ by
some number of Dehn twists, 
and those may be cancelled by rigidly rotating $F_1$.

We have thus constructed
a regular homotopy between $i$ and $i\circ h$ such that the 
number mod 2 of quadruple points,
is the sum of the number occurring in the 
$\bar{F_k}$s in case $h$ is orientation preserving, 
and the sum plus 1, in case
$h$ is orientation reversing.
In other words $Q(i,i\circ h)=Q(i_1,i_1\circ h_1)+Q(i_2,i_2\circ h_2)+\ep(h)$.
By the induction hypothesis, $Q(i_k,i_k\circ h_k)=\hp(h_k)$, $k=1,2$. 
Let $n_k={\mathrm{genus}}(F_k)$
and notice $n=n_1+n_2$, $\ep(h_k)=\ep(h)$ and
$\r(h_*-Id)=\r({h_1}_*-Id)+\r({h_2}_*-Id)$ and so
$\p(h_*)=\p({h_1}_*) + \p({h_2}_*)$. So finally:
$Q(i,i\circ h)=\hp(h_1)+\hp(h_2)+\ep(h)=
\p({h_1}_*)+(n_1+1)\ep(h)+\p({h_2}_*)+(n_2+1)\ep(h)+\ep(h)=
\p(h_*)+(n_1+n_2+1)\ep(h)  = \hp(h)$.

We now deal with the special $U$-map generator which appears
when $\g=2,\Ai=0$. Since $\Ai =0$, then by 
Theorem \pr{Pi1}(2,3) and Lemma \pr{bb} as before, we may assume
$i$ is an embedding whose image is two embedded tori connected to each
other with a tube, and such that a half revolution around some line in $\E$
maps $i(F)$ onto itself and interchanges the two tori.
Let $h:F\to F$ be the map such that $i\circ h$ is the final embedding
of this half revolution,
then $h$ is orientation preserving and so $h\in\mi$.
Take a circle $c$ in one of the tori
with $g^i([c])=1$, then $h(c)$ lies in the other torus, 
and so $[c]\neq h_*([c])$ and $[c]\cdot h_*([c])=0$.
It follows that $[c]$ and $h_*([c])$ are not in the same $V_k$ of 
the definition of $U$-map, so $h_*$ must be a $U$-map, and so
$h$ is indeed a $U$-map on $F$.
Since $h^2=Id$ then by Lemma \pr{pum}, $\p(h_*)=0$ 
and so $\hp(h)=\p(h_*)+\ep(h)=0$.
Since there is a rigid 
rotation between $i$ and $i\circ h$, $Q(i,i\circ h)=0$.
This completes the proof of Theorem \pr{main1}.


\begin{thebibliography}{StoPC}

\bibitem[C]{C}
C.C. Chevalley: ``The Algebraic Theory of Spinors.'' 
Columbia University Press 1954. Also reprinted in:
C. Chevalley: ``The Algebraic Theory of Spinors and Clifford Algebras, 
Collected Works, Vol. 2'' 
Springer-Verlag 1997.

\bibitem[H]{H}
M.W. Hirsch: ``Immersions of manifolds.'' 
\emph{Trans. Amer. Math. Soc.} 93 (1959) 242--276. 

\bibitem[L]{L}
W.B.R. Lickorish: ``An introduction to knot theory.'' 
Graduate Texts in Mathematics 175. Springer-Verlag, New York, 1997.

\bibitem[MB]{MB}
N. Max, T. Banchoff: ``Every Sphere Eversion Has a Quadruple Point.''
\emph{Contributions to Analysis and Geometry}, 
John Hopkins University Press,
(1981), 191-209. 

\bibitem[N]{N}
T. Nowik:
``Quadruple points of regular homotopies of surfaces in 3-manifolds.''
\emph{Topology} - to appear.
(may be viewed at: http://www.math.columbia.edu/$\tilde{\ }$tahl/quadruple.ps)

\bibitem[P]{P}
U. Pinkall: ``Regular homotopy classes of immersed surfaces.''
\emph{Topology} 24 (1985) No.4, 421--434.

\bibitem[S]{S}
S. Smale: ``A classification of immersions of the two-sphere.'' 
\emph{Trans. Amer. Math. Soc.} 90 (1958) 281-290. 

\end{thebibliography}
\end{document}